\documentclass[11pt,oneside,english]{amsart}
\usepackage[T1]{fontenc}
\usepackage[latin9]{inputenc}
\usepackage{babel}
\usepackage{amstext}
\usepackage{amsthm}
\usepackage{amssymb}
\usepackage{geometry}
\usepackage[bookmarks=true,bookmarksnumbered=false,bookmarksopen=false,
 breaklinks=false,pdfborder={0 0 1},backref=false,colorlinks=false]
 {hyperref}

\makeatletter
\numberwithin{equation}{section}
\numberwithin{figure}{section}
\theoremstyle{plain}
\newtheorem{thm}{\protect\theoremname}[section]
\theoremstyle{plain}
\newtheorem{lem}[thm]{\protect\lemmaname}
\theoremstyle{plain}
\newtheorem{prop}[thm]{\protect\propositionname}

\IfFileExists{lmodern.sty}{\usepackage{lmodern}}{}

\makeatother

\providecommand{\lemmaname}{Lemma}
\providecommand{\propositionname}{Proposition}
\providecommand{\theoremname}{Theorem}

\begin{document}
\title{The pair correlation of $\alpha n+\beta\log n$ modulo one}
\author{Nadav Yesha}
\address{Department of Mathematics, University of Haifa, 3103301 Haifa, Israel}
\begin{abstract}
We show that for any non-Liouville irrational $\alpha$ and for any
real $\beta\ne0$, the sequence $\left(\alpha n+\beta\log n\right)_{n=1}^{\infty}$
has Poissonian pair correlation modulo one.
\end{abstract}

\maketitle

\section{Introduction}

The pair correlation function of a sequence $\left(x_{n}\right)_{n=1}^{\infty}$
of real numbers, defined for any $s>0$ by
\begin{equation}
R_{N}\left(s\right)=\frac{1}{N}\#\left\{ 1\le m\ne n\le N:\,\left\Vert x_{m}-x_{n}\right\Vert \le\frac{s}{N}\right\} \label{eq:PC_def}
\end{equation}
where $\left\Vert x\right\Vert $ denotes the distance from $x$ to
the nearest integer, is an important local statistic which helps measure
the pseudorandomness of the sequence $\left(x_{n}\right)_{n=1}^{\infty}$
modulo one. The sequence is said to have Poissonian pair correlation
if for any $s>0$,
\begin{equation}
\lim_{N\to\infty}R_{N}\left(s\right)=2s.\label{eq:Poisson_limit}
\end{equation}

The goal of this paper is to prove Poissonian pair correlation for
the sequence
\[
x_{n}=\alpha n+\beta\log n
\]
for any non-Liouville irrational $\alpha$, and any $\beta\ne0$.
\begin{thm}
\label{thm:MainThm}Let $\alpha$ be a non-Liouville irrational, and
let $\beta\ne0$ be real. Then the sequence $x_{n}=\alpha n+\beta\log n$
has Poissonian pair correlation, i.e., for any $s>0$, we have
\[
\lim_{N\to\infty}R_{N}\left(s\right)=2s.
\]
\end{thm}

Taking $\beta=\alpha$, Theorem \ref{thm:MainThm} shows that $x_{n}=\alpha\left(n+\log n\right)$
has Poissonian pair correlation for every non-Liouville irrational
$\alpha$, providing a deterministic version of the recent metric
result of Tan, Zhang and Zhou \cite{TZZ1}.

Since replacing $\left(\alpha,\beta\right)$ by $\left(-\alpha,-\beta\right)$
leaves $R_{N}\left(s\right)$ unchanged, it suffices to consider $\beta>0$,
which we will assume throughout.

\subsection{Background}

Recall that the sequence $\left(x_{n}\right)_{n=1}^{\infty}$ is said
to be uniformly distributed modulo one (u.d. mod 1) if for any $0\le a<b\le1$
we have
\[
\lim_{N\to\infty}\frac{1}{N}\#\left\{ 1\le n\le N:\,\left\{ x_{n}\right\} \in[a,b)\right\} =b-a,
\]
where $\left\{ x\right\} $ denotes the fractional part of $x$. This
classical notion dates back to the work of Weyl \cite{Weyl} and has
been extensively studied for more than a century. There are many examples
of sequences which are known to be u.d. mod 1. Most famously, Weyl
showed \cite{Weyl} that every sequence of the form $x_{n}=p\left(n\right)$
where $p\left(x\right)=a_{d}x^{d}+\dots+a_{1}x+a_{0}$ is u.d. mod
1 if at least one of the non-constant coefficients $a_{1},\dots,a_{d}$
is irrational. In particular, the Kronecker sequence $x_{n}=\alpha n$
is uniformly distributed modulo one for any irrational $\alpha$.
Further examples include sequences of the form $x_{n}=\alpha n^{\theta}$
where $\alpha\ne0$ and $\theta>0$ is not an integer. There are also
results which hold in a metric sense, i.e., for (Lebesgue) almost
every parameter in some parameterised family of sequences. A remarkable
metric theorem is also due to Weyl, who showed \cite{Weyl} that for
any strictly increasing sequence of integers $\left(a_{n}\right)_{n=1}^{\infty}$,
the sequence $x_{n}=\alpha a_{n}$ is u.d. mod $1$ for almost all
$\alpha\in\mathbb{R}$.

The notion of u.d. mod 1 only gives a very coarse measure of the pseudorandomness
of a sequence. A better way to understand if a sequence has random-like
behaviour is to consider its local statistics at the scale $1/N$
(the mean gap between its first $N$ elements modulo one). These local
statistics are governed by the correlations of the sequence, the most
accessible of them being the pair correlation (\ref{eq:PC_def}).
In the random model where the sequence $\left(x_{n}\right)_{n=1}^{\infty}$
consists of uniform i.i.d. random variables in the unit interval,
the pair correlation function indeed satisfies the limit (\ref{eq:Poisson_limit})
almost surely. Thus, Poissonian pair correlation is an indication
of pseudorandomness which holds at the local scale. It is strictly
stronger than uniform distribution modulo one (see \cite{ALP,GrepstadLarcher}).

The bulk of the results concerning Poissonian pair correlation have
been obtained in a metric sense. Rudnick and Sarnak proved \cite{RudnickSarnak}
Poissonian pair correlation for the polynomial sequences $x_{n}=\alpha n^{d}$
($d\ge2$ integer) for almost all $\alpha$. More recently, for each
fixed $\theta>1$, Aistleitner, El-Baz and Munsch \cite{AEM} proved
Poissonian pair correlation for the sequence $x_{n}=\alpha n^{\theta}$
for almost all $\alpha$. Rudnick and Technau \cite{RudnickTechnau}
established the corresponding result for $0<\theta<1$. For each fixed
$\alpha>0$, Shubin \cite{Shubin} proved Poissonian pair correlation
for the sequence $x_{n}=\alpha n^{\theta}$ for almost all $\theta\in\left(0,\frac{3}{5}\right)\cup\left(3,\infty\right)$,
extending earlier work of Technau and Yesha \cite{TechnauYesha}.
Aistleitner, Baker, Technau and Yesha proved \cite{ABTY} that for
almost all $\alpha>1$, the sequence $x_{n}=\alpha^{n}$ has Poissonian
correlations of every order $k\ge2$. Very recently, Tan, Zhang and
Zhou obtained Poissonian pair correlation for almost all $\alpha\in\mathbb{R}$
for the following two sequences: $x_{n}=\alpha\left(n+\log n\right)$
\cite{TZZ1} and $x_{n}=\alpha n\log n$ \cite{TZZ2}, answering two
open questions of Aistleitner, El-Baz and Munsch from \cite{AEM}.

Establishing Poissonian pair correlation for natural sequences with
prescribed parameters (i.e., not in the metric sense) is substantially
harder, and so far has only been done for a few examples. For the
quadratic polynomials $x_{n}=\alpha n^{2}$, Heath-Brown constructed
\cite{HeathBrown} a dense set of coefficients $\alpha$ such that
Poissonian pair correlation holds, although Poissonian pair correlation
remains unproved for simple prescribed coefficients such as $\alpha=\sqrt{2}$.
El-Baz, Marklof and Vinogradov \cite{EMV} showed Poissonian pair
correlation for the sequence $x_{n}=n^{1/2}$ with perfect-square
indices removed (this is perhaps somewhat surprising, as this sequence
has a non-Poissonian gap distribution, as was proved by Elkies and
McMullen \cite{ElkiesMcMullen}). More recently, Lutsko, Sourmelidis
and Technau showed \cite{LST} Poissonian pair correlation for the
slowly growing sequences $x_{n}=\alpha n^{\theta}$ with $\alpha>0$
and $0<\theta<\frac{14}{41}=\frac{1}{3}+0.008\dots$, which was improved
by Radziwi\l\l{} and Shubin \cite{RadziwillShubin} to $0<\theta<\frac{43}{117}=\frac{1}{3}+0.034\dots$.
For even more slowly growing sequences of the form $x_{n}=\alpha\left(\log n\right)^{A}$
where $\alpha>0$ and $A>1$, the correlations of all orders were
shown to be Poissonian by Lutsko and Technau \cite{LutskoTechnau}.
A more complicated recent example is due to Hauke \cite{Hauke}, who
constructed certain sequences $\left(a_{n}\right)_{n=1}^{\infty}$
of integers arising from the set of rough numbers such that for any
badly approximable $\alpha$, the correlations of all orders of the
sequence $x_{n}=\alpha a_{n}$ are Poissonian.

In our paper, we deal with the sequence $x_{n}=\alpha n+\beta\log n$.
It is easy to show that for any irrational $\alpha$ this sequence
is u.d. mod 1 (see \cite{Ohkubo,TichyTurnwald} for bounds for the
discrepancy of these sequences for non-Liouville $\alpha$). One can
think of this sequence as a logarithmic perturbation of the Kronecker
sequence. While the three-gap theorem, together with \cite[Theorem 1]{LarcherStockinger},
implies that the pair correlation of the Kronecker sequence is not
Poissonian, our result shows that for a non-Liouville irrational $\alpha$,
such logarithmic perturbations introduce enough pseudorandomness to
restore Poissonian pair correlations. While similar questions can
also be studied in a more general setting, we chose to focus on logarithmic
perturbations for which the inverse relation used in the proof has
an explicit form.

\subsection{Outline of the argument and organisation of the paper}

The paper is organised as follows. In Section \ref{sec:AuxLemmas},
we obtain the auxiliary estimates used in the proof. Section \ref{sec:Equid}
is devoted to proving the equidistribution result in Proposition \ref{prop:Eq_Prop},
which is key in our argument. By symmetry, we restrict to $m>n$ and
write $m=n+h$, so that our indices lie in the triangle $n\ge1,$
$h\ge1$ and $n+h\le N$. In the bulk of the triangle (by which we
mean away from the axes), the pair correlation condition means that
for a suitable integer $k$ with $k-\alpha h>0$, $n$ lies in an
interval around the root of the equation
\[
\alpha h+\beta\log\left(1+\frac{h}{x_{0}}\right)=k,
\]
which is
\[
x_{0}=\frac{h}{e^{\frac{k-\alpha h}{\beta}}-1}.
\]
We study these roots modulo one, jointly with the parameters appearing
in the count. After a Fourier expansion, Diophantine approximation
allows us to split the resulting sums into arithmetic progressions.
Choosing the progression step and the derivative order appropriately,
we obtain cancellation from van der Corput's derivative estimates.
Section \ref{sec:Bulk} uses this equidistribution result to evaluate
the pair correlation count in the bulk of the triangle. In Section
\ref{sec:axes}, we show that the contribution from small $n$ or
small $h$ is negligible. Finally, the above estimates are combined
in Section \ref{sec:Proof} where we conclude the proof of Theorem
\ref{thm:MainThm}.

\section{\label{sec:AuxLemmas}Auxiliary lemmas}

In this section, we prove a few auxiliary lemmas. We will say that
an irrational $\alpha$ is of Diophantine type $\tau\ge2$ if there
exists a constant $c=c\left(\alpha\right)>0$ such that for any $p\in\mathbb{Z}$
and $q\in\mathbb{N}$ we have
\begin{equation}
\left|\alpha-\frac{p}{q}\right|\ge\frac{c}{q^{\tau}}.\label{eq:Type_tau}
\end{equation}
In particular, any non-Liouville irrational $\alpha\in\mathbb{R}$
is of finite type $\tau$ for some $\tau\ge2$ (note that there are
other conventions for defining the type in the literature, see e.g.
\cite{KuipersNiederreiter} for another common convention).

The first lemma is concerned with Diophantine approximation of irrationals
of finite type.
\begin{lem}
\label{lem:Dioph}Let $\alpha\in\mathbb{R}$ be irrational of finite
type $\tau\ge2$. For any $Q\ge1$, there exist integers $p$ and
$1\le q\ll_{\alpha}Q^{\tau-1}$ such that
\[
\frac{1}{Q}\le\left|\alpha q-p\right|<\frac{2}{Q}.
\]
\end{lem}

\begin{proof}
By Dirichlet's approximation theorem, there exist integers $p_{0},q_{0}$
with $1\le q_{0}\le Q$ such that
\[
\left|\alpha q_{0}-p_{0}\right|<\frac{1}{Q}.
\]
Let $m=\left\lceil \frac{1}{Q\left|\alpha q_{0}-p_{0}\right|}\right\rceil $.
Then
\[
\frac{1}{Q}\le m\left|\alpha q_{0}-p_{0}\right|<\frac{1}{Q}+\left|\alpha q_{0}-p_{0}\right|<\frac{2}{Q}.
\]
Thus, if we let $p=mp_{0},$ $q=mq_{0}$, then
\[
\frac{1}{Q}\le\left|\alpha q-p\right|<\frac{2}{Q}.
\]
Moreover, since $\left|\alpha q_{0}-p_{0}\right|\ge\frac{c}{q_{0}^{\tau-1}},$
we have
\[
q<q_{0}+\frac{q_{0}}{Q\left|\alpha q_{0}-p_{0}\right|}\le Q+\frac{q_{0}^{\tau}}{cQ}\ll_{\alpha}Q^{\tau-1}.
\]
This proves our claim.
\end{proof}
Note that the integers $p$ and $q$ in Lemma \ref{lem:Dioph} need
not be coprime, which is harmless for our purpose.

The next lemma concerns the distribution of the fractional parts of
the Kronecker sequence in short intervals.
\begin{lem}
\label{lem:DiscrepancyLemma}Let $\alpha\in\mathbb{R}$ be irrational
of finite type $\tau\ge2$, and let $J\subseteq\mathbb{R}/\mathbb{Z}$
be an interval. For any integers $M,H\ge1$, and for any $\epsilon>0$,
we have
\[
\#\left\{ M<n\le M+H:\,\left\{ \alpha n\right\} \in J\right\} =\left|J\right|H+O_{\alpha,\epsilon}\left(H^{1-\frac{1}{\tau-1}+\epsilon}\right)
\]
where the implied constant in the error term depends only on $\alpha$
and $\epsilon$.
\begin{proof}
Let
\begin{equation}
D_{N}=\frac{1}{N}\sup_{I\subseteq\mathbb{R}/\mathbb{Z}}\Bigl|\#\left\{ 1\le n\le N:\,\left\{ \alpha n\right\} \in I\right\} -\left|I\right|N\Bigr|\label{eq:Discrepancy}
\end{equation}
denote the discrepancy of the sequence $\left(\alpha n\right)_{n=1}^{\infty}$,
where the supremum is taken over all half-open intervals in $\mathbb{R}/\mathbb{Z}$.
Clearly, the definition (\ref{eq:Discrepancy}) is equivalent to the
definition where the supremum is taken over half open intervals of
the form $[a,b)$ where $0\le a<b\le1$ (see, e.g., \cite[p. 88]{KuipersNiederreiter}).
Note that $D_{N}$ is invariant under translations of the sequence,
and therefore
\begin{align*}
\Bigl|\#\left\{ M<n\le M+H:\,\left\{ \alpha n\right\} \in J\right\} -\left|J\right|H\Bigr| & =\Bigl|\#\left\{ 1\le n\le H:\,\left\{ \alpha n+\alpha M\right\} \in J\right\} -\left|J\right|H\Bigr|\\
 & \le HD_{H}.
\end{align*}
Since $\alpha$ is of finite type $\tau$, by standard discrepancy
bounds (see, e.g. \cite[p. 123, Theorem 3.2]{KuipersNiederreiter},
taking into consideration their alternative convention for the type
of $\alpha$) we have $D_{H}=O_{\alpha,\epsilon}\left(H^{-\frac{1}{\tau-1}+\epsilon}\right)$
for any $\epsilon>0$, which yields our claim.
\end{proof}
\end{lem}

We will also need the following version of van der Corput's $k$th
derivative test:
\begin{lem}[{\cite[Theorems 5.9, 5.11, 5.13]{Titchmarsh}}]
\label{lem:VdC}Let $\phi\in C^{k}\left(\left[a,b\right]\right)$
be a real valued function, where $k\ge2$. Assume that there exist
$\lambda_{k}>0$ and $C\ge1$ such that
\[
\lambda_{k}\le\left|\phi^{\left(k\right)}\left(x\right)\right|\le C\lambda_{k},
\]
for all $x\in\left[a,b\right]$. Denote $K=2^{k-1}$ and $M=\max\left\{ b-a,1\right\} $.
Then
\begin{equation}
\sum_{a<n\le b}e^{2\pi i\phi\left(n\right)}\ll C^{2/K}M\lambda_{k}^{\frac{1}{2K-2}}+M^{1-\frac{2}{K}}\lambda_{k}^{-\frac{1}{2K-2}}\label{eq:osc_sum}
\end{equation}
where the implied constant is absolute.
\end{lem}

Note that in \cite{Titchmarsh}, the result is stated for $b-a\ge1.$
If $b-a<1$, then the sum in (\ref{eq:osc_sum}) has at most one term,
and hence the bound holds trivially, since $\lambda_{k}^{\frac{1}{2K-2}}+\lambda_{k}^{-\frac{1}{2K-2}}\ge2$.

Applying partial summation, we immediately deduce the following smooth
weight version of Lemma \ref{lem:VdC}.
\begin{lem}
\label{lem:VdCWeight}Let $w\in C^{1}\left(\left[a,b\right]\right)$
be complex valued, and assume that $\phi$ satisfies the conditions
of Lemma \ref{lem:VdC}. Then
\[
\sum_{a<n\le b}w\left(n\right)e^{2\pi i\phi\left(n\right)}\ll\left(\sup_{x\in\left[a,b\right]}\left|w\left(x\right)\right|+\int_{a}^{b}\left|w'\left(x\right)\right|dx\right)\left(C^{2/K}M\lambda_{k}^{\frac{1}{2K-2}}+M^{1-\frac{2}{K}}\lambda_{k}^{-\frac{1}{2K-2}}\right)
\]
where the implied constant is absolute.
\end{lem}

\section{\label{sec:Equid}An equidistribution result}

From now on, we fix a non-Liouville irrational $\alpha$ of finite
type $\tau$, and $\beta>0$. For $t>0$, we denote
\begin{equation}
F\left(t\right)=\frac{1}{e^{\frac{t}{\beta}}-1}.\label{eq:F_def}
\end{equation}
The goal of this section is to prove an equidistribution result, which
will be a key ingredient of our proof.
\begin{prop}
\label{prop:Eq_Prop}Let $g\in C_{c}^{\infty}\left(\left(0,\infty\right)^{2}\times\mathbb{R}/\mathbb{Z}\right)$
be real valued. We have
\begin{equation}
\lim_{N\to\infty}\frac{1}{N}\sum_{\substack{k\in\mathbb{Z},\,h\ge1,\\
k-\alpha h>0
}
}g\left(k-\alpha h,\frac{h}{N},hF\left(k-\alpha h\right)\right)=\int_{0}^{\infty}\int_{0}^{\infty}\int_{0}^{1}g\left(u,t,x\right)\,dxdtdu.\label{eq:eq_sum}
\end{equation}
\end{prop}

In order to prove Proposition \ref{prop:Eq_Prop}, we will expand
$g$ in the third variable into its Fourier series. For any nonzero
frequency $l$, the required cancellation will follow from the following
key lemma:
\begin{lem}
\label{lem:Exponential_lemma}Let $\psi\in C_{c}^{\infty}\left(\left(0,\infty\right)^{2}\right)$
be complex valued, and let $l\in\mathbb{Z}\setminus\left\{ 0\right\} .$
As $N\to\infty$, we have
\begin{equation}
\sum_{\substack{k\in\mathbb{Z},\,h\ge1,\\
k-\alpha h>0
}
}\psi\left(k-\alpha h,\frac{h}{N}\right)e^{2\pi ilhF\left(k-\alpha h\right)}=O_{l,\alpha,\beta,\psi}\left(N^{1-\sigma}\right)\label{eq:exp_sum}
\end{equation}
where $\sigma=\frac{1}{\lfloor\tau+1\rfloor\left(2^{\lfloor\tau+1\rfloor+1}-2\right)}$.
\end{lem}

\begin{proof}
Let $r=\lfloor\tau+1\rfloor+1$ and let $Q=N^{\frac{1}{r-1}}$. Note
that $Q^{\tau}=o\left(N\right)$. By Lemma \ref{lem:Dioph}, there
exist integers $p,q$ with $q\ll Q^{\tau-1}$ so that
\begin{equation}
\left|\alpha q-p\right|\asymp\frac{1}{Q}.\label{eq:Dirichlet_approx}
\end{equation}
Denote $\varepsilon=q\alpha-p$, so that by (\ref{eq:Dirichlet_approx})
we have $\left|\varepsilon\right|\asymp\frac{1}{Q}.$ We decompose
the summation over $h$ in (\ref{eq:exp_sum}) into arithmetic progressions:
\begin{align}
h & =a+jq,\hspace{1em}0\le a<q.\label{eq:arith_prog}
\end{align}
Denote
\[
k=c+jp,
\]
so that
\begin{equation}
k-\alpha h=c-\alpha a-j\varepsilon.\label{eq:k-alpha_h}
\end{equation}

Since $h\asymp N$, (\ref{eq:arith_prog}) implies that for any fixed
$a$, we have $j=O\left(1+\frac{N}{q}\right)$. Since $k-\alpha h$
is bounded, (\ref{eq:k-alpha_h}) implies that the number of possible
values of $c$ is $O\left(1+\left|\varepsilon\right|\frac{N}{q}\right)$.
Hence, the number of possible choices of $\left(a,c\right)$ is
\[
O\left(q+\left|\varepsilon\right|N\right)\ll Q^{\tau-1}+\frac{N}{Q}\ll\frac{N}{Q}.
\]

Let $\left(s_{-},s_{+}\right)\times\left(t_{-},t_{+}\right)$ with
$s_{-},t_{-}>0$ be a rectangle which contains the support of $\psi$.
For any fixed $a$ and $c$, let
\[
J=J_{a,c}=\left\{ x\in\mathbb{R}:\,s_{-}\le c-\alpha a-x\varepsilon\le s_{+},\,t_{-}\le\frac{a+xq}{N}\le t_{+}\right\}
\]
which is an interval (possibly empty) of length $O\left(\frac{1}{\left|\varepsilon\right|}\right)=O\left(Q\right)$.
The indices $j$ which contribute to (\ref{eq:exp_sum}) clearly belong
to $J$. Denote this part of the sum (\ref{eq:exp_sum}) by
\[
S_{J}=\sum_{j\in J}\psi\left(c-\alpha a-j\varepsilon,\frac{a+jq}{N}\right)e^{2\pi il\left(a+jq\right)F\left(c-\alpha a-j\varepsilon\right)}.
\]
The phase function is then
\[
\phi\left(x\right)=l\left(a+qx\right)F\left(c-\alpha a-\varepsilon x\right),
\]
and its $r$th derivative is clearly
\begin{equation}
\phi^{\left(r\right)}\left(x\right)=l\left(a+qx\right)\left(-\varepsilon\right)^{r}F^{\left(r\right)}\left(c-\alpha a-\varepsilon x\right)+lqr\left(-\varepsilon\right)^{r-1}F^{\left(r-1\right)}\left(c-\alpha a-\varepsilon x\right).\label{eq:phi_r}
\end{equation}
Note that for each $x\in J$, $c-\alpha a-\varepsilon x$ lies in
a compact interval in $\left(0,\infty\right)$, and $a+qx\asymp N$.
Also note that for every integer $j\ge0$ and every  $t>0$, we have
\[
\left(-1\right)^{j}F^{\left(j\right)}\left(t\right)=\sum_{m=1}^{\infty}\left(\frac{m}{\beta}\right)^{j}e^{-\frac{mt}{\beta}}>0,
\]
and hence by continuity $F^{\left(j\right)}\left(t\right)$ is bounded
away from zero on compact intervals of $\left(0,\infty\right)$. We
therefore see that
\[
\left|\frac{lqr\left(-\varepsilon\right)^{r-1}F^{\left(r-1\right)}\left(c-\alpha a-\varepsilon x\right)}{l\left(a+qx\right)\left(-\varepsilon\right)^{r}F^{\left(r\right)}\left(c-\alpha a-\varepsilon x\right)}\right|\ll\frac{q}{\left|\varepsilon\right|N}\ll\frac{Q^{\tau}}{N}=o\left(1\right)
\]
uniformly over the contributing $\left(a,c\right)$ and $x\in J_{a,c}$,
so that the first term in (\ref{eq:phi_r}) dominates the second term,
and thus
\[
\left|\phi^{\left(r\right)}\left(x\right)\right|\asymp N\left|\varepsilon\right|^{r}\asymp\frac{N}{Q^{r}}=\frac{1}{Q}.
\]

By Lemma \ref{lem:VdCWeight} with the weight $w\left(x\right)=\psi\left(c-\alpha a-\varepsilon x,\frac{a}{N}+\frac{q}{N}x\right),$
for which
\[
\sup\limits_{x\in J}\left|w\left(x\right)\right|+\int_{J}\left|w'\left(x\right)\right|dx\ll1+Q\left(\left|\varepsilon\right|+\frac{q}{N}\right)\ll1,
\]
we deduce that for $K=2^{r-1}$, we have
\[
S_{J}\ll Q\left(\frac{1}{Q}\right)^{\frac{1}{2K-2}}+Q^{1-\frac{2}{K}}\left(\frac{1}{Q}\right)^{-\frac{1}{2K-2}}.
\]
Summing over the $O\left(\frac{N}{Q}\right)$ options for $\left(a,c\right)$,
we conclude that
\begin{align*}
\sum_{\substack{k\in\mathbb{Z},\,h\ge1,\\
k-\alpha h>0
}
}\psi\left(k-\alpha h,\frac{h}{N}\right)e^{2\pi ilhF\left(k-\alpha h\right)} & \ll NQ^{-\frac{1}{2K-2}}+NQ^{-\frac{2}{K}+\frac{1}{2K-2}}\ll NQ^{-\frac{1}{2K-2}}
\end{align*}
since
\[
\frac{2}{2K-2}<\frac{2}{K}.
\]
Since $NQ^{-\frac{1}{2K-2}}=N^{1-\sigma}$ with $\sigma=\frac{1}{\left(r-1\right)\left(2K-2\right)}$
this completes the proof.
\end{proof}
We now prove Proposition \ref{prop:Eq_Prop}.
\begin{proof}[Proof of Proposition \ref{prop:Eq_Prop}]
Fix $u$ and $t$, and expand $g$ in the $x$ variable into its
Fourier series, so that
\[
g\left(u,t,x\right)=\sum_{l\in\mathbb{Z}}\psi_{l}\left(u,t\right)e^{2\pi ilx},
\]
where
\[
\psi_{l}\left(u,t\right)=\int_{0}^{1}g\left(u,t,s\right)e^{-2\pi ils}\,ds.
\]
By the smoothness of $g$, the coefficients $\psi_{l}\left(u,t\right)$
decay rapidly in $l$, so that as $T\to\infty$, we have
\[
g\left(u,t,x\right)=\sum_{\left|l\right|\le T}\psi_{l}\left(u,t\right)e^{2\pi ilx}+R_{T}\left(u,t,x\right),
\]
where $R_{T}\left(u,t,x\right)=O\left(T^{-1}\right)$ uniformly in
$u,t$ and $x$, and vanishes whenever $\left(u,t\right)$ lies outside
a fixed compact set in $\left(0,\infty\right)^{2}$, independent of
$T$. The sum in (\ref{eq:eq_sum}) is supported on $h\asymp N$,
and for each $h$ there is only a bounded number of possible values
of $k$ contributing to the sum. Thus,
\begin{align}
\frac{1}{N}\sum_{\substack{k\in\mathbb{Z},\,h\ge1,\\
k-\alpha h>0
}
}g\left(k-\alpha h,\frac{h}{N},hF\left(k-\alpha h\right)\right) & =\frac{1}{N}\sum_{\left|l\right|\le T}\sum_{\substack{k\in\mathbb{Z},\,h\ge1,\\
k-\alpha h>0
}
}\psi_{l}\left(k-\alpha h,\frac{h}{N}\right)e^{2\pi ilhF\left(k-\alpha h\right)}\label{eq:g_Fourier}\\
 & +O\left(T^{-1}\right).
\end{align}
For any $l\ne0,$ Lemma \ref{lem:Exponential_lemma} gives
\[
\sum_{\substack{k\in\mathbb{Z},\,h\ge1,\\
k-\alpha h>0
}
}\psi_{l}\left(k-\alpha h,\frac{h}{N}\right)e^{2\pi ilhF\left(k-\alpha h\right)}=O_{l}\left(N^{1-\sigma}\right),
\]
so that
\begin{equation}
\frac{1}{N}\sum_{0<\left|l\right|\le T}\sum_{\substack{k\in\mathbb{Z},\,h\ge1,\\
k-\alpha h>0
}
}\psi_{l}\left(k-\alpha h,\frac{h}{N}\right)e^{2\pi ilhF\left(k-\alpha h\right)}=O_{T}\left(N^{-\sigma}\right).\label{eq:nonzero_l}
\end{equation}
For $l=0$, extend $\psi_{0}\left(\cdot,t\right)$ by zero to $\mathbb{R}$.
Poisson summation in the $k$ variable yields
\[
\frac{1}{N}\sum_{\substack{k\in\mathbb{Z},\,h\ge1,\\
k-\alpha h>0
}
}\psi_{0}\left(k-\alpha h,\frac{h}{N}\right)=\frac{1}{N}\sum_{j\in\mathbb{Z}}\sum_{h\ge1}e^{-2\pi ijh\alpha}a_{j}\left(h\right)
\]
where
\begin{align*}
a_{j}\left(h\right) & =\int_{0}^{\infty}\psi_{0}\left(u,\frac{h}{N}\right)e^{-2\pi iju}\,du.
\end{align*}
As before, by the rapid decay of $a_{j}\left(h\right)$ in $j$, we
have
\begin{align}
\frac{1}{N}\sum_{\substack{k\in\mathbb{Z},\,h\ge1,\\
k-\alpha h>0
}
}\psi_{0}\left(k-\alpha h,\frac{h}{N}\right) & =\frac{1}{N}\sum_{\left|j\right|\le T}\sum_{h\ge1}e^{-2\pi ijh\alpha}a_{j}\left(h\right)+O\left(T^{-1}\right).\label{eq:psi_Fourier}
\end{align}
For any $j\ne0$, we use partial summation. Write
\[
A_{j}\left(t\right)=\int_{0}^{\infty}\psi_{0}\left(u,t\right)e^{-2\pi iju}\,du,
\]
then since for any $H_{1}<H_{2}$ we have
\[
\left|\sum_{H_{1}<h\le H_{2}}e^{-2\pi ijh\alpha}\right|\ll\frac{1}{\left\Vert j\alpha\right\Vert },
\]
we get that
\begin{equation}
\sum_{h\ge1}e^{-2\pi ijh\alpha}a_{j}\left(h\right)=\sum_{h\ge1}e^{-2\pi ijh\alpha}A_{j}\left(\frac{h}{N}\right)\ll\frac{1}{\left\Vert j\alpha\right\Vert }\left(\left\Vert A_{j}\right\Vert _{\infty}+\left\Vert A_{j}'\right\Vert _{1}\right)\ll_{j}1.\label{eq:nonzero_j}
\end{equation}
The term $j=0$ is a Riemann sum, which contributes
\begin{equation}
\frac{1}{N}\sum_{h\ge1}\int_{0}^{\infty}\psi_{0}\left(u,\frac{h}{N}\right)\,du=\int_{0}^{\infty}\int_{0}^{\infty}\int_{0}^{1}g\left(u,t,x\right)\,dxdtdu+O\left(N^{-1}\right).\label{eq:main_term}
\end{equation}
Thus, collecting (\ref{eq:g_Fourier}), (\ref{eq:nonzero_l}), (\ref{eq:psi_Fourier}),
(\ref{eq:nonzero_j}), and (\ref{eq:main_term}), we conclude that
\begin{align*}
\frac{1}{N}\sum_{\substack{k\in\mathbb{Z},\,h\ge1,\\
k-\alpha h>0
}
}g\left(k-\alpha h,\frac{h}{N},hF\left(k-\alpha h\right)\right) & =\int_{0}^{\infty}\int_{0}^{\infty}\int_{0}^{1}g\left(u,t,x\right)\,dxdtdu\\
 & +O_{T}\left(N^{-\sigma}\right)+O\left(T^{-1}\right).
\end{align*}
Thus, there exists a constant $C>0$ such that for any $T>0$ we have
\begin{align*}
\int_{0}^{\infty}\int_{0}^{\infty}\int_{0}^{1}g\left(u,t,x\right)\,dxdtdu-\frac{C}{T} & \le\liminf_{N\to\infty}\frac{1}{N}\sum_{\substack{k\in\mathbb{Z},\,h\ge1,\\
k-\alpha h>0
}
}g\left(k-\alpha h,\frac{h}{N},hF\left(k-\alpha h\right)\right)\\
 & \le\limsup_{N\to\infty}\frac{1}{N}\sum_{\substack{k\in\mathbb{Z},\,h\ge1,\\
k-\alpha h>0
}
}g\left(k-\alpha h,\frac{h}{N},hF\left(k-\alpha h\right)\right)\\
 & \le\int_{0}^{\infty}\int_{0}^{\infty}\int_{0}^{1}g\left(u,t,x\right)\,dxdtdu+\frac{C}{T}
\end{align*}
so that taking $T\to\infty$ proves our claim.
\end{proof}

\section{\label{sec:Bulk}The correlation in the bulk}

We will now see how the equidistribution result we proved in Proposition
\ref{prop:Eq_Prop} connects with the pair correlation problem. We
assume for now that $m>n$, and write $m=n+h$ with $1\le n\le N$,
$h\ge1$, and $1\le h+n\le N$, so that
\[
x_{n+h}-x_{n}=\alpha h+\beta\log\left(1+\frac{h}{n}\right).
\]
We thus need to estimate the number of pairs $\left(n,h\right)$ in
the triangle
\[
T_{N}=\left\{ \left(x,y\right):\,1\le x\le N,\,y\ge1,\,1\le x+y\le N\right\}
\]
such that
\[
\left\Vert \alpha h+\beta\log\left(1+\frac{h}{n}\right)\right\Vert \le\frac{s}{N}.
\]
We will first estimate a smooth count of the number of such elements.
Let $W\in C_{c}^{\infty}\left(\left(0,\infty\right)^{2}\right)$ and
$f\in C_{c}^{\infty}\left(\mathbb{R}\right)$ be real valued, and
denote the smooth counting function by
\begin{equation}
S_{N}=S_{W,f,N}=\frac{1}{N}\sum_{n,h\ge1}W\left(\frac{n}{N},\frac{h}{N}\right)\sum_{k\in\mathbb{Z}}f\left(N\left(\alpha h-k+\beta\log\left(1+\frac{h}{n}\right)\right)\right).\label{eq:S_N}
\end{equation}
Our main goal will be to prove a Poissonian limit for $S_{N}$:
\begin{lem}
\label{prop:S_N_Poisson}As $N\to\infty$, we have
\[
\lim_{N\to\infty}S_{N}=\int_{0}^{\infty}\int_{0}^{\infty}W\left(x,t\right)\,dxdt\int_{-\infty}^{\infty}f\left(z\right)\,dz.
\]
\end{lem}

Note that since $W$ is compactly supported in $\left(0,\infty\right)^{2}$,
as we range over $n$ and $h$, $\log\left(1+\frac{h}{n}\right)$
is contained in a compact interval inside $\left(0,\infty\right)$.
Thus, since
\begin{equation}
k-\alpha h=\beta\log\left(1+\frac{h}{n}\right)+O\left(\frac{1}{N}\right),\label{eq:phi(n)_small}
\end{equation}
$k-\alpha h$ also lies in a compact interval in $\left(0,\infty\right)$
for sufficiently large $N$. Hence, for any $h$, there is only a
bounded number of possible values of $k$ contributing to the sum
(\ref{eq:S_N})$.$ Define $g\in C_{c}^{\infty}\left(\left(0,\infty\right)^{2}\times\mathbb{R}/\mathbb{Z}\right)$
by
\begin{equation}
g\left(u,t,x\right)=W\left(tF\left(u\right),t\right)\sum_{n\in\mathbb{Z}}f\left(\frac{1}{tF'\left(u\right)}\left(n-x\right)\right).\label{eq:g_def}
\end{equation}
To see that $g$ is smooth and compactly supported in $\left(0,\infty\right)^{2}\times\mathbb{R}/\mathbb{Z}$,
note that by the compact supports of $W$ and $f$ and by the definition
of $F$, the function $W\left(tF\left(u\right),t\right)$ is compactly
supported in $\left(0,\infty\right)^{2}$ and the sum over $n$ is
locally finite.
\begin{lem}
\label{lem:S_N_series}As $N\to\infty$, we have
\begin{equation}
S_{N}=\frac{1}{N}\sum_{\substack{k\in\mathbb{Z},\,h\ge1,\\
k-\alpha h>0
}
}g\left(k-\alpha h,\frac{h}{N},hF\left(k-\alpha h\right)\right)+O_{\alpha,\beta,W,f}\left(\frac{1}{N}\right).\label{eq:S_N_expansion}
\end{equation}
\end{lem}

\begin{proof}
We compare the summands over the union of their supports. As discussed
above, in the sum defining (\ref{eq:S_N}) we can restrict to $n,h\asymp N$
and $k-\alpha h$ lying in a fixed compact interval in $\left(0,\infty\right)$.
Note that the same is true in the sum involved in $(\ref{eq:S_N_expansion})$
since $g$ has compact support in $\left(0,\infty\right)^{2}$ with
respect to its first two variables and since $f$ is compactly supported.
We can therefore assume that $n,h,k$ are variables satisfying these
conditions. Let $x_{0}=hF\left(k-\alpha h\right),$ so that also $x_{0}\asymp N$.
Denote
\[
\Phi\left(x\right)=\alpha h-k+\beta\log\left(1+\frac{h}{x}\right),
\]
and note that $\Phi\left(x_{0}\right)=0$ and that
\begin{equation}
x=hF\left(\Phi\left(x\right)+k-\alpha h\right).\label{eq:inverse_formula}
\end{equation}
Hence,
\begin{equation}
n-x_{0}=hF\left(\Phi\left(n\right)+k-\alpha h\right)-hF\left(k-\alpha h\right)=O\left(N\left|\Phi\left(n\right)\right|\right).\label{eq:n-x_0}
\end{equation}
By (\ref{eq:phi(n)_small}), for $n,h,k$ involved in the sum (\ref{eq:S_N})
we have $\Phi\left(n\right)=O\left(\frac{1}{N}\right)$, so that $(\ref{eq:n-x_0})$
gives
\begin{equation}
n-x_{0}=O\left(1\right).\label{eq:n-x_0_bound}
\end{equation}
For $n,h,k$ involved in the sum (\ref{eq:S_N_expansion}), (\ref{eq:n-x_0_bound})
also clearly follows from the compact support of $f$, so we may assume
it holds. In particular, we have
\begin{equation}
\frac{n}{N}=\frac{x_{0}}{N}+O\left(\frac{1}{N}\right).\label{eq:n/N}
\end{equation}

Differentiating (\ref{eq:inverse_formula}), we get
\[
\Phi'\left(x\right)=\frac{1}{hF'\left(\Phi\left(x\right)+k-\alpha h\right)}
\]
so that $\Phi'\left(x_{0}\right)=\frac{1}{hF'\left(k-\alpha h\right)}$.
We also have the formulas
\[
\Phi'\left(x\right)=-\beta\frac{h}{x\left(x+h\right)},\hspace{1em}\Phi''\left(x\right)=\frac{\beta h\left(2x+h\right)}{x^{2}\left(x+h\right)^{2}}
\]
so that$\left|\Phi'\left(x\right)\right|\asymp\frac{1}{N}$ and $\left|\Phi''\left(x\right)\right|\asymp\frac{1}{N^{2}}$
for any $x\asymp N$. Taylor expanding $\Phi$ around $x_{0}$, evaluating
at $x=n$, and using (\ref{eq:n-x_0_bound}) to bound the error term,
we get
\begin{equation}
\Phi\left(n\right)=\frac{n-x_{0}}{hF'\left(k-\alpha h\right)}+O\left(\frac{1}{N^{2}}\right).\label{eq:g_eq}
\end{equation}

The number of allowed triples in the above discussions is $O\left(N\right)$,
since there are $O\left(N\right)$ choices for $h$, and given $h$
only $O\left(1\right)$ choices for $k$, and by (\ref{eq:n-x_0_bound})
only $O\left(1\right)$ choices for $n$. Since the derivatives of
$W$ and $f$ are bounded, the estimates (\ref{eq:g_eq}) and (\ref{eq:n/N})
imply
\begin{equation}
W\left(\frac{n}{N},\frac{h}{N}\right)f\left(N\Phi\left(n\right)\right)-W\left(\frac{x_{0}}{N},\frac{h}{N}\right)f\left(\frac{N\left(n-x_{0}\right)}{hF'\left(k-\alpha h\right)}\right)=O\left(\frac{1}{N}\right)\label{eq:replacement}
\end{equation}
uniformly over these triples.

Summing (\ref{eq:replacement}) over the $O\left(N\right)$ contributing
triples, we obtain
\[
S_{N}=\frac{1}{N}\sum_{\substack{k\in\mathbb{Z},\,h\ge1,\\
k-\alpha h>0
}
}W\left(\frac{x_{0}}{N},\frac{h}{N}\right)\sum_{n\in\mathbb{Z}}f\left(\frac{N}{hF'\left(k-\alpha h\right)}\left(n-x_{0}\right)\right)+O\left(\frac{1}{N}\right).
\]
Recalling the definition (\ref{eq:g_def}) of the function $g$, the
claimed result follows.
\end{proof}
We are now ready to prove Lemma \ref{prop:S_N_Poisson}.
\begin{proof}[Proof of Lemma \ref{prop:S_N_Poisson}]
By Lemma \ref{lem:S_N_series} and Proposition \ref{prop:Eq_Prop},
we have
\begin{align*}
\lim_{N\to\infty}S_{N} & =\lim_{N\to\infty}\frac{1}{N}\sum_{\substack{k\in\mathbb{Z},\,h\ge1,\\
k-\alpha h>0
}
}g\left(k-\alpha h,\frac{h}{N},hF\left(k-\alpha h\right)\right)\\
 & =\int_{0}^{\infty}\int_{0}^{\infty}\int_{0}^{1}g\left(u,t,x\right)\,dxdtdu\\
 & =\int_{0}^{\infty}\int_{0}^{\infty}W\left(tF\left(u\right),t\right)\int_{-\infty}^{\infty}f\left(\frac{-x}{tF'\left(u\right)}\right)\,dxdtdu\\
 & =\int_{0}^{\infty}\int_{0}^{\infty}-tF'\left(u\right)W\left(tF\left(u\right),t\right)\,dtdu\int_{-\infty}^{\infty}f\left(z\right)\,dz\\
 & =\int_{0}^{\infty}\int_{0}^{\infty}W\left(x,t\right)\,dxdt\int_{-\infty}^{\infty}f\left(z\right)\,dz
\end{align*}
as claimed.
\end{proof}
Let $\eta_{1},\eta_{2}>0,$ such that $\eta_{1}+\eta_{2}<1$, and
let
\[
T_{\eta_{1},\eta_{2}}=\left\{ \left(x,y\right):\,\eta_{1}\le x\le1,\eta_{2}\le y,\,x+y\le1\right\} ,
\]
which is a triangle of area $\frac{1}{2}\left(1-\eta_{1}-\eta_{2}\right)^{2}$,
and let
\[
S_{N,\eta_{1},\eta_{2}}=\frac{1}{N}\#\left\{ \left(n,h\right)\in\mathbb{N}^{2}:\,\left(\frac{n}{N},\frac{h}{N}\right)\in T_{\eta_{1},\eta_{2}},\,\left\Vert x_{n+h}-x_{n}\right\Vert \le\frac{s}{N}\right\}
\]
be the non-smooth count of the pair correlation in the bulk of the
triangle $T_{N}$.
\begin{lem}
\label{lem:Bulk_lemma}For any $\eta_{1},\eta_{2}>0$ such that $\eta_{1}+\eta_{2}<1$,
we have
\[
\lim_{N\to\infty}S_{N,\eta_{1},\eta_{2}}=\left(1-\eta_{1}-\eta_{2}\right)^{2}s.
\]
\end{lem}

\begin{proof}
Let $\epsilon>0$. Since both sets $T_{\eta_{1},\eta_{2}}$ and $\left[-s,s\right]$
have boundary of Lebesgue measure zero, we can approximate their indicator
functions from below and from above by non-negative, smooth functions
$W_{\pm}\left(x,t\right)$ and $f_{\pm}\left(z\right)$, where $W_{\pm}$
are compactly supported in $\left(0,\infty\right)^{2}$ and $f_{\pm}$
are compactly supported in $\mathbb{R}$, so that
\[
0\le W_{-}\le1_{T_{\eta_{1},\eta_{2}}}\le W_{+},\hspace{1em}0\le f_{-}\le1_{\left[-s,s\right]}\le f_{+},
\]
and
\[
\left|\int_{0}^{\infty}\int_{0}^{\infty}W_{\pm}\left(x,t\right)\,dxdt\int_{-\infty}^{\infty}f_{\pm}\left(z\right)\,dz-\left(1-\eta_{1}-\eta_{2}\right)^{2}s\right|<\epsilon.
\]
For $N>2s$ we have
\[
S_{W_{-},f_{-},N}\le S_{N,\eta_{1},\eta_{2}}\le S_{W_{+},f_{+},N}.
\]
Hence, by Lemma \ref{prop:S_N_Poisson}, we have
\[
\left(1-\eta_{1}-\eta_{2}\right)^{2}s-\epsilon\le\liminf_{N\to\infty}S_{N,\eta_{1},\eta_{2}}\le\limsup_{N\to\infty}S_{N,\eta_{1},\eta_{2}}\le\left(1-\eta_{1}-\eta_{2}\right)^{2}s+\epsilon
\]
which proves our claim since $\epsilon>0$ is arbitrarily small.
\end{proof}

\section{\label{sec:axes}Small indices}

We now estimate how pairs $\left(n,h\right)$ where either $n$ or
$h$ is small contribute to the pair correlation, and show that their
contribution is negligible. We will first deal with the case when
$n$ is large and $h$ may be small.
\begin{lem}
\label{lem:small_h}Let $\eta_{1}>0$, $c>0$, and $1\le H_{0}\le N$.
For any $0<\epsilon<\frac{1}{\tau\left(\tau-1\right)}$, we have
\begin{gather}
\#\left\{ \left(n,h\right)\in\mathbb{N}^{2}:\,\eta_{1}N\le n\le cN,\,1\le h\le H_{0},\,\left\Vert x_{n+h}-x_{n}\right\Vert \le\frac{s}{N}\right\} \label{eq:n_h_count}\\
=O_{s,\eta_{1},c,\alpha,\beta,\epsilon}\left(H_{0}+N^{1-\frac{1}{\tau\left(\tau-1\right)}+\epsilon}\right).\nonumber
\end{gather}
\end{lem}

\begin{proof}
Let $1\le H\le H_{0}$, and assume first that $H\le h<2H$. If $\left\Vert x_{n+h}-x_{n}\right\Vert =\left\Vert \alpha h+\beta\log\left(1+\frac{h}{n}\right)\right\Vert \le\frac{s}{N}$
then
\begin{equation}
\left\Vert \alpha h\right\Vert \le\beta\log\left(1+\frac{h}{n}\right)+\frac{s}{N}\ll\frac{H}{N}.\label{eq:alpha_a_small}
\end{equation}
By (\ref{eq:Type_tau}), we have $H\ge CN^{\frac{1}{\tau}}$ for some
constant $C>0.$ Moreover, by Lemma \ref{lem:DiscrepancyLemma}, the
number of $H\le h<2H$ such that (\ref{eq:alpha_a_small}) holds is
$O\left(\frac{H^{2}}{N}+H^{1-\frac{1}{\tau-1}+\epsilon}\right)$ for
any $\epsilon>0$. For any such fixed $h,$ we have to count
\begin{equation}
\sum_{k\in\mathbb{Z}}\#\left\{ \eta_{1}N\le n\le cN:\,\left|\Phi\left(n\right)\right|\le\frac{s}{N}\right\} \label{eq:n_sum_fixed_h}
\end{equation}
where $\Phi\left(x\right)=\alpha h-k+\beta\log\left(1+\frac{h}{x}\right)$.
Note that there are $O\left(1\right)$ possible values of $k$ for
each fixed $h$, since $\log\left(1+\frac{h}{n}\right)=O\left(1\right)$.
For each $k$, the derivative of $\Phi$ is
\[
\left|\Phi'\left(x\right)\right|=\beta\frac{h}{x\left(x+h\right)}\asymp\frac{H}{N^{2}}
\]
for $x\asymp N$, and therefore the number of $n$ contributing to
the sum (\ref{eq:n_sum_fixed_h}) is $O\left(1+\frac{N}{H}\right)=O\left(\frac{N}{H}\right).$
In total, the contribution of $\left(n,h\right)$ for $H\le h<2H$
to the count (\ref{eq:n_h_count}) is thus
\[
O\left(\frac{N}{H}\left(\frac{H^{2}}{N}+H^{1-\frac{1}{\tau-1}+\epsilon}\right)\right)=O\left(H+NH^{-\frac{1}{\tau-1}+\epsilon}\right).
\]

Now divide the interval $\left[CN^{\frac{1}{\tau}},H_{0}\right]$
into dyadic intervals. The total contribution of all $\left(n,h\right)$
is then bounded by
\[
\sum_{CN^{\frac{1}{\tau}}\le2^{k}\le H_{0}}\left(2^{k}+N2^{k\left(-\frac{1}{\tau-1}+\epsilon\right)}\right)\ll H_{0}+N^{1+\frac{1}{\tau}\left(-\frac{1}{\tau-1}+\epsilon\right)}
\]
which gives our claim.
\end{proof}
The following lemma will deal with small values of $n$.
\begin{lem}
\label{lem:Small_n_lemma}Let $1\le H\le N$. We have
\[
\#\left\{ \left(n,m\right)\in\mathbb{N}^{2}:\,1\le n\le H,\,1\le m\le N,\,\left\Vert x_{m}-x_{n}\right\Vert \le\frac{s}{N}\right\} =O_{\alpha,\beta,s}\left(\left(HN\right)^{1/2}\right).
\]
\end{lem}

\begin{proof}
For $0\le j<N$, denote $I_{j}=[\frac{j}{N},\frac{j+1}{N}).$ Given
a finite index set $T$, let $v_{T}\in\mathbb{R}^{N}$ denote the
vector with coordinates
\[
v_{T}^{j}=\#\left\{ n\in T:\,\left\{ x_{n}\right\} \in I_{j}\right\} =\sum_{n\in T}1_{I_{j}}\left(\left\{ x_{n}\right\} \right).
\]
Note that
\begin{align*}
\left\Vert v_{T}\right\Vert _{2}^{2} & =\sum_{j=0}^{N-1}\sum_{n,m\in T}1_{I_{j}}\left(\left\{ x_{m}\right\} \right)1_{I_{j}}\left(\left\{ x_{n}\right\} \right)\\
 & =\#T+2\#\left\{ \left(n,m\right)\in T^{2}:\,n<m,\,\exists_{0\le j<N},\left\{ x_{m}\right\} ,\left\{ x_{n}\right\} \in I_{j}\right\} .
\end{align*}
Let $1\le M\le N$ and $T\subseteq\left\{ M,M+1,\dots,2M-1\right\} $,
and note that if $n,m\in T$ and $\left\{ x_{m}\right\} ,\left\{ x_{n}\right\} $
belong to the same bin $I_{j}$, then
\[
\left\Vert x_{m}-x_{n}\right\Vert \le\frac{1}{N}\le\frac{1}{M}.
\]
Thus, Lemma \ref{lem:small_h} with $N=M$, $H_{0}=M,$ $\eta_{1}=1$,
$c=2$ and $s=1$ yields
\begin{equation}
\left\Vert v_{T}\right\Vert _{2}^{2}\le M+2\#\left\{ \left(n,h\right)\in\mathbb{N}^{2}:\,M\le n\le2M,\,1\le h\le M,\,\left\Vert x_{n+h}-x_{n}\right\Vert \le\frac{1}{M}\right\} \ll M.\label{eq:l2_norm_v}
\end{equation}
By decomposing the indices $\left\{ 1,\dots,H\right\} $ into dyadic
blocks $T_{k}$, the triangle inequality and the bound (\ref{eq:l2_norm_v})
applied for each block give
\[
\left\Vert v_{\left\{ 1,\dots,H\right\} }\right\Vert _{2}\le\sum_{2^{k}\le H}\left\Vert v_{T_{k}}\right\Vert _{2}\ll\sum_{2^{k}\le H}\left(2^{k}\right)^{1/2}\ll H^{1/2}.
\]
 Likewise, $\left\Vert v_{\left\{ 1,\dots,N\right\} }\right\Vert _{2}\ll N^{1/2}.$
Now, if $\left\{ x_{n}\right\} \in I_{j}$, $\left\{ x_{m}\right\} \in I_{k}$,
and $\left\Vert x_{m}-x_{n}\right\Vert \le\frac{s}{N}$, then clearly
$k=j+r\,\left(\text{mod}\,N\right)$ where $\left|r\right|\le\lceil s\rceil+1$.
Hence, by Cauchy-Schwarz,
\begin{gather*}
\#\left\{ \left(n,m\right)\in\mathbb{N}^{2}:\,1\le n\le H,\,1\le m\le N,\,\left\Vert x_{m}-x_{n}\right\Vert \le\frac{s}{N}\right\} \\
\le\sum_{\left|r\right|\le\lceil s\rceil+1}\sum_{j=0}^{N-1}v_{\left\{ 1,\dots,H\right\} }^{j}v_{\left\{ 1,\dots,N\right\} }^{j+r\,\left(\text{mod}\,N\right)}\ll\left\Vert v_{\left\{ 1,\dots,H\right\} }\right\Vert _{2}\left\Vert v_{\left\{ 1,\dots,N\right\} }\right\Vert _{2}\ll\left(HN\right)^{1/2}
\end{gather*}
as claimed.
\end{proof}

\section{\label{sec:Proof}Proof of Theorem \ref{thm:MainThm}}

We now have all that we need to complete the proof of Theorem \ref{thm:MainThm}.
\begin{proof}[Proof of Theorem \ref{thm:MainThm}]
Let $\eta_{1},\eta_{2}>0$ such that $\eta_{1}+\eta_{2}<1$. By Lemma
\ref{lem:Bulk_lemma} with $h=m-n$, we have
\[
\lim_{N\to\infty}\frac{1}{N}\#\left\{ \eta_{1}N\le n<m\le N:\,m-n\ge\eta_{2}N,\,\left\Vert x_{m}-x_{n}\right\Vert \le\frac{s}{N}\right\} =\left(1-\eta_{1}-\eta_{2}\right)^{2}s.
\]
By Lemma \ref{lem:small_h} again with $h=m-n$, we have
\[
\limsup_{N\to\infty}\frac{1}{N}\#\left\{ \eta_{1}N\le n<m\le N:\,m-n<\eta_{2}N,\,\left\Vert x_{m}-x_{n}\right\Vert \le\frac{s}{N}\right\} \le C_{\eta_{1}}\cdot\eta_{2}
\]
where $C_{\eta_{1}}$ is a constant depending on $\eta_{1}$ (but
not on $\eta_{2}$). Finally, by Lemma \ref{lem:Small_n_lemma}, we
have
\[
\limsup_{N\to\infty}\frac{1}{N}\#\left\{ 1\le n<m\le N:\,n<\eta_{1}N,\,\left\Vert x_{m}-x_{n}\right\Vert \le\frac{s}{N}\right\} \le C\eta_{1}^{1/2}
\]
where $C$ is a constant (independent of $\eta_{1}$ or $\eta_{2}$).
Thus,
\begin{gather*}
2\left(1-\eta_{1}-\eta_{2}\right)^{2}s\le\liminf_{N\to\infty}R_{N}\left(s\right)\le\limsup_{N\to\infty}R_{N}\left(s\right)\le2\left(1-\eta_{1}-\eta_{2}\right)^{2}s+2C_{\eta_{1}}\eta_{2}+2C\eta_{1}^{1/2}.
\end{gather*}
If we first fix $\eta_{1}>0$ and let $\eta_{2}\to0$, we conclude
that
\[
2\left(1-\eta_{1}\right)^{2}s\le\liminf_{N\to\infty}R_{N}\left(s\right)\le\limsup_{N\to\infty}R_{N}\left(s\right)\le2\left(1-\eta_{1}\right)^{2}s+2C\eta_{1}^{1/2}.
\]
We now let $\eta_{1}\to0$ and conclude that
\[
\lim_{N\to\infty}R_{N}\left(s\right)=2s
\]
which concludes the proof of Theorem \ref{thm:MainThm}.
\end{proof}

\subsection*{Acknowledgement of AI use}

The arguments in this paper were developed through discussions with
GPT-6 Astra (OpenAI), starting from the author's initial ideas. The
author wrote the manuscript and independently verified all its contents.
The model was also used for proofreading and review of the manuscript.
The author takes full responsibility for the contents of the paper.

\end{document}